\documentclass{amsart}
\usepackage[square,numbers]{natbib}
\newtheorem{thm}{Theorem}
\newtheorem{lem}[thm]{Lemma}
\newtheorem{cor}[thm]{Corollary}

\theoremstyle{definition}

\newtheorem{ex}[thm]{Example}

\providecommand{\abs}[1]{\lvert#1\rvert}

\providecommand{\norm}[1]{\lVert#1\rVert}

\begin{document}
\title[Duality]{Some duality results for equivalence couplings and total variation}

\author{Luca Pratelli}
\address{Luca Pratelli, Accademia Navale, viale Italia 72, 57100 Livorno,
Italy} \email{pratel@mail.dm.unipi.it}

\author{Pietro Rigo}
\address{Pietro Rigo, Dipartimento di Scienze Statistiche ``P. Fortunati'', Universit\`a di Bologna, via delle Belle Arti 41, 40126 Bologna, Italy}
\email{pietro.rigo@unibo.it}



\begin{abstract}
Let $(\Omega,\mathcal{F})$ be a standard Borel space and $\mathcal{P}(\mathcal{F})$ the collection of all probability measures on $\mathcal{F}$. Let $E\subset\Omega\times\Omega$ be a measurable equivalence relation, that is, $E\in\mathcal{F}\otimes\mathcal{F}$ and the relation on $\Omega$ defined as $x\sim y$ $\Leftrightarrow$ $(x,y)\in E$ is reflexive, symmetric and transitive. It is shown that there are two $\sigma$-fields $\mathcal{G}_0$ and $\mathcal{G}_1$ on $\Omega$ such that, for all $\mu,\,\nu\in\mathcal{P}(\mathcal{F})$,
$$\inf_{P\in\Gamma(\mu,\nu)}(1-P(E))=\norm{\mu-\nu}_{\mathcal{G}_1}\quad\text{and}\quad\min_{P\in\Gamma(\mu,\nu_0)}(1-P(E))=\norm{\mu-\nu}_{\mathcal{G}_0}.$$
Here, $\nu_0\in\mathcal{P}(\mathcal{F})$ is a suitable probability measure satisfying $\nu_0=\nu$ on $\mathcal{G}_0$. Moreover, $\mathcal{G}_0\subset\mathcal{F}$ while $\mathcal{G}_1\subset\widehat{\mathcal{F}}$, where $\widehat{\mathcal{F}}$ is the universally measurable $\sigma$-field with respect to $\mathcal{F}$. However, for all $\mu,\,\nu\in\mathcal{P}(\mathcal{F})$, there is a $\sigma$-field $\mathcal{G}(\mu,\nu)\subset\mathcal{F}$ such that
$$\inf_{P\in\Gamma(\mu,\nu)}(1-P(E))=\norm{\mu-\nu}_{\mathcal{G}(\mu,\nu)}.$$
\end{abstract}

\maketitle

\section{Introduction}\label{intro}
\noindent Throughout, $(\Omega,\mathcal{F})$ is a measurable space, $\mathcal{P}(\mathcal{F})$ the collection of all probability measures on $\mathcal{F}$, and $E\subset\Omega\times\Omega$ a measurable equivalence relation. This means that $E\in\mathcal{F}\otimes\mathcal{F}$ and the relation on $\Omega$ defined as
$$x\sim y\quad\Leftrightarrow\quad (x,y)\in E$$
is reflexive, symmetric and transitive.

\medskip

\noindent The following notion of duality has been recently introduced by Jaffe \cite{JAFFE}. Given a sub-$\sigma$-field $\mathcal{G}\subset\mathcal{F}$, the pair $(E,\mathcal{G})$ is said to satisfy {\em strong duality} if
$$\min_{P\in\Gamma(\mu,\nu)}(1-P(E))=\norm{\mu-\nu}_\mathcal{G}\quad\quad\text{for all }\mu,\,\nu\in\mathcal{P}(\mathcal{F}).$$
Here, as usual, $\Gamma(\mu,\nu)$ is the set of all probability measures on $\mathcal{F}\otimes\mathcal{F}$ with marginals $\mu$ and $\nu$ and the notation ``$\min$" asserts that the infimum is actually achieved. Moreover,
$$\norm{\mu-\nu}_\mathcal{G}=\sup_{A\in\mathcal{G}}\,\abs{\mu(A)-\nu(A)}$$
is the total variation between $\mu$ and $\nu$ on $\mathcal{G}$.

\medskip

\noindent In addition to be intriguing from the foundational point of view, strong duality is useful in some probabilistic frameworks. Examples concern stochastic calculus, point processes, and random sequence simulation; see Section 2 of \cite{JAFFE}.

\medskip

\noindent Say that $E$ is {\em strongly dualizable} if $(E,\mathcal{G})$ satisfies strong duality for {\em some} sub-$\sigma$-field $\mathcal{G}\subset\mathcal{F}$. Various conditions for $E$ to be strongly dualizable are given in \cite{JAFFE} but no measurable equivalence relation which fails to be strongly dualizable is known to date. This suggests the {\em conjecture} that, under mild conditions on $(\Omega,\mathcal{F})$ (say $(\Omega,\mathcal{F})$ is a standard Borel space), every measurable equivalence relation is strongly dualizable.

\medskip

\noindent This note focus on strong duality and includes three results. Let
$$\mathcal{G}_0=\bigl\{A\in\mathcal{F}:\,1_A(x)=1_A(y)\text{ for all }(x,y)\in E\bigr\}.$$
Such a $\mathcal{G}_0$ is a sub-$\sigma$-field of $\mathcal{F}$ which plays a special role as regards strong duality. In fact, $E$ is strongly dualizable if and only if $(E,\mathcal{G}_0)$ satisfies strong duality; see \cite[Prop. 3.15]{JAFFE}. Our first result is that, for all $\mu,\,\nu\in\mathcal{P}(\mathcal{F})$, there is $\nu_0\in\mathcal{P}(\mathcal{F})$ satisfying
\begin{gather*}
\nu_0=\nu\text{ on }\mathcal{G}_0\quad\text{and}\quad\min_{P\in\Gamma(\mu,\nu_0)}(1-P(E))=\norm{\mu-\nu}_{\mathcal{G}_0}.
\end{gather*}
Roughly speaking, the above condition means that strong duality is always true up to changing one between $\mu$ and $\nu$ out of $\mathcal{G}_0$. This is quite reasonable, in a sense, for $\norm{\mu-\nu}_{\mathcal{G}_0}$ only involves the restrictions of $\mu$ and $\nu$ on $\mathcal{G}_0$.

\medskip

\noindent Next, suppose $(\Omega,\mathcal{F})$ is a standard Borel space and denote by $\widehat{\mathcal{F}}$ the collection of those subsets of $\Omega$ which are universally measurable with respect to $\mathcal{F}$; see Section \ref{prel}. Define also
$$\mathcal{G}_1=\bigl\{A\in\widehat{\mathcal{F}}:\,1_A(x)=1_A(y)\text{ for all }(x,y)\in E\bigr\}.$$
This time, $\mathcal{G}_1$ is not a sub-$\sigma$-field of $\mathcal{F}$. However, by our second result, one obtains
\begin{gather}\label{m8a}
\inf_{P\in\Gamma(\mu,\nu)}(1-P(E))=\norm{\mu-\nu}_{\mathcal{G}_1}\quad\quad\text{for all }\mu,\,\nu\in\mathcal{P}(\mathcal{F}).
\end{gather}
In addition, the $\inf$ is achieved if $P$ is allowed to be finitely additive. Precisely,
\begin{gather*}
\min_{P\in M(\mu,\nu)}(1-P(E))=\norm{\mu-\nu}_{\mathcal{G}_1}\quad\quad\text{for all }\mu,\,\nu\in\mathcal{P}(\mathcal{F})
\end{gather*}
where $M(\mu,\nu)$ is the collection of finitely additive probabilities on $\mathcal{F}\otimes\mathcal{F}$ with marginals $\mu$ and $\nu$.

\medskip

\noindent Finally, for each $B\subset\Omega$, define
$$\mathcal{G}_B=\bigl\{A\in\mathcal{F}:1_A(x)=1_A(y)\text{ for all }(x,y)\in E\cap (B\times B)\bigr\}.$$
Then, for all $\mu,\,\nu\in\mathcal{P}(\mathcal{F})$, there is a set $B\in\mathcal{F}$ such that
\begin{gather*}
\mu(B)=\nu(B)=1\quad\text{and}\quad\inf_{P\in\Gamma(\mu,\nu)}(1-P(E))=\norm{\mu-\nu}_{\mathcal{G}_B}.
\end{gather*}
If compared with \eqref{m8a}, the latter result has the advantage that $\mathcal{G}_B$ is a sub-$\sigma$-field of $\mathcal{F}$ but the disadvantage that $\mathcal{G}_B$ is not universal, for it depends on the pair $(\mu,\nu)$. Note also that $\mu(B)=\nu(B)=1$ and $E\cap (B\times B)$ is a measurable equivalence relation on $B$. Therefore, for fixed $(\mu,\nu)$, one can replace $\Omega$ with $B$ and $E$ with $E\cap (B\times B)$. After doing this, everything works as regards the total variation side of strong duality.

\section{Preliminaries}\label{prel}
\noindent In this section, we introduce some further notation and recall a few known facts.

\medskip

\noindent Let $(S,\mathcal{E})$ be a measurable space. Then, $\mathcal{P}(\mathcal{E})$ denotes the set of probability measures on $\mathcal{E}$ and $M(\mathcal{E})$ the set of bounded $\mathcal{E}$-measurable functions $f:S\rightarrow\mathbb{R}$. We write
$$\mu(f)=\int f\,d\mu\quad\quad\text{whenever }\mu\in\mathcal{P}(\mathcal{E})\text{ and }f\in M(\mathcal{E}).$$
Moreover, we let
$$\widehat{\mathcal{E}}=\bigcap_{\mu\in\mathcal{P}(\mathcal{E})}\overline{\mathcal{E}}^\mu$$
where $\overline{\mathcal{E}}^\mu$ is the completion of $\mathcal{E}$ with respect to $\mu$. The elements of $\widehat{\mathcal{E}}$ are usually called {\em universally measurable} with respect to $\mathcal{E}$. With a slight abuse of notation, for each $\mu\in\mathcal{P}(\mathcal{E})$, the unique extension of $\mu$ to $\widehat{\mathcal{E}}$ is still denoted by $\mu$.

\medskip

\noindent If $T$ is any topological space, $\mathcal{B}(T)$ denotes the Borel $\sigma$-field. We say that $T$ is {\em Polish} if its topology is induced by a distance $d$ such that $(T,d)$ is a complete separable metric space. Moreover, the  measurable space $(S,\mathcal{E})$ is a {\em standard Borel space} if $\mathcal{E}=\mathcal{B}(S)$ for some Polish topology on $S$.

\medskip

\noindent A probability $\mu\in\mathcal{P}(\mathcal{E})$ is {\em perfect} if, for any $\mathcal{E}$-measurable function $f:S\rightarrow\mathbb{R}$, there is a Borel set $B\in\mathcal{B}(\mathbb{R})$ such that $B\subset f(S)$ and $\mu(f\in B)=1$. In a sense, perfectness is a non-topological analogous of the notion of tightness. In fact, if $S$ is separable metric and $\mathcal{E}=\mathcal{B}(S)$, then $\mu$ is perfect if and only if it is tight. In particular, each element of $\mathcal{P}(\mathcal{E})$ is perfect whenever $(S,\mathcal{E})$ is a standard Borel space. We refer to \cite{RAM1979} for more on perfect probability measures.

\medskip

\noindent As regards duality theory in mass transportation, we just mention a result by Ramachandran and Ruschendorf \cite{RAMRUS1995}. For more information, the interested reader is referred to \cite{BLS}, \cite{BNT17}, \cite{KELL}, \cite{RIG20} and references therein. Given $\mu,\,\nu\in\mathcal{P}(\mathcal{E})$, let $\Gamma(\mu,\nu)$ be the collection of probability measures $P$ on $\mathcal{E}\otimes\mathcal{E}$ with marginals $\mu$ and $\nu$, i.e.
$$P(A\times S)=\mu(A)\quad\text{and}\quad P(S\times A)=\nu(A)\quad\text{for all }A\in\mathcal{E}.$$
Moreover, let $c:S\times S\rightarrow\mathbb{R}$ be a bounded measurable cost function. (Boundedness of $c$ is generally superfluous and has been assumed for the sake of simplicity only). A {\em primal minimizer}, or an {\em optimal coupling}, is a probability measure $P\in\Gamma(\mu,\nu)$ such that $P(c)\le Q(c)$ for each $Q\in\Gamma(\mu,\nu)$. For a primal minimizer to exist, it suffices that $S$ is separable metric, $\mathcal{E}=\mathcal{B}(S)$, $\mu$ and $\nu$ are perfect, and the cost $c$ is lower semi-continuous. To state the duality result, we denote by $L$ the set of pairs $(f,g)$ satisfying
$$f,\,g\in M(\mathcal{E})\quad\text{and}\quad f(x)+g(y)\le c(x,y)\text{ for all }(x,y)\in S\times S.$$
Then, in view of \cite{RAMRUS1995}, one obtains
$$\inf_{P\in\Gamma(\mu,\nu)}P(c)=\sup_{(f,g)\in L}\bigl\{\mu(f)+\nu(g)\bigr\}$$
provided at least one between $\mu$ and $\nu$ is perfect.

\medskip

\noindent Finally, let $\mathcal{D}\subset\mathcal{E}$ be a sub-$\sigma$-field and $\mu,\,\nu\in\mathcal{P}(\mathcal{E})$. For any measure $\gamma$ on $\mathcal{E}$, we write $\gamma|\mathcal{D}$ to denote the restriction of $\gamma$ on $\mathcal{D}$. The total variation between $\mu$ and $\nu$ on $\mathcal{D}$ is
$$\norm{\mu-\nu}_\mathcal{D}=\sup_{A\in\mathcal{D}}\,\abs{\mu(A)-\nu(A)}=\sup_{
                                                                               \begin{array}{c}
                                                                                 f\in M(\mathcal{D}) \\
                                                                                 0\le f\le 1 \\
                                                                               \end{array}}\,\abs{\mu(f)-\nu(f)}.$$
For our purposes, two remarks are in order. First, $\norm{\cdot}_\mathcal{D}$ can be written as
\begin{gather}\label{g8m4e3}
\norm{\mu-\nu}_\mathcal{D}=\mu(A)-\nu(A)\quad\quad\text{for a suitable }A\in\mathcal{D}.
\end{gather}
Define in fact
\begin{gather}\label{f67y9k}
\lambda=\mu+\nu,\quad f=\frac{d\,(\mu|\mathcal{D})}{d\,(\lambda|\mathcal{D})}\quad\text{and}\quad g=\frac{d\,(\nu|\mathcal{D})}{d\,(\lambda|\mathcal{D})}.
\end{gather}
Letting $A=\bigl\{f>g\bigr\}$, one obtains $A\in\mathcal{D}$ (since $f$ and $g$ are $\mathcal{D}$-measurable) and
\begin{gather*}
\norm{\mu-\nu}_\mathcal{D}=\int\left(f-g\right)^+\,d\lambda=\int_A\left(f-g\right)\,d\lambda=\mu(A)-\nu(A).
\end{gather*}
The second remark consists in the following lemma, which slightly improves some well known facts; see e.g. \cite[Prop. 3.1]{BPR10} and \cite[Lem. 2.1]{SET}.

\begin{lem}\label{kl8}
Let $\mathcal{D}\subset\mathcal{E}$ be a sub-$\sigma$-field and $\mu,\,\nu\in\mathcal{P}(\mathcal{E})$. Then, there are a probability space $(\Phi,\mathcal{A},\mathbb{P})$ and two measurable maps
$X,\,Y:(\Phi,\mathcal{A})\rightarrow (S,\mathcal{E})$ such that
\begin{gather*}
\mathbb{P}(X\in A)=\mu(A)\text{ for all }A\in\mathcal{E},\quad\mathbb{P}(Y\in A)=\nu(A)\text{ for all }A\in\mathcal{D},
\\\bigl\{X\ne Y\bigr\}\in\mathcal{A}\quad\text{and}\quad\mathbb{P}(X\ne Y)=\norm{\mu-\nu}_\mathcal{D}.
\end{gather*}
\end{lem}

\begin{proof} Suppose first $\mu|\mathcal{D}=\nu|\mathcal{D}$. Let $\Phi=S\times S$, $\mathcal{C}=\mathcal{E}\otimes\mathcal{E}$ and $X(a,b)=a$ and $Y(a,b)=b$ for all $(a,b)\in S\times S$. Define also
$$\mathbb{Q}(C)=\mu\bigl\{x\in S:(x,x)\in C\bigr\}\quad\quad\text{for all }C\in\mathcal{C}.$$
Then, $\mathbb{Q}(X\in A)=\mathbb{Q}(Y\in A)=\mu(A)$ for all $A\in\mathcal{E}$. In particular, since $\mu=\nu$ on $\mathcal{D}$, one obtains $\mathbb{Q}(Y\in A)=\nu(A)$ for all $A\in\mathcal{D}$. Moreover, under $\mathbb{Q}$, the inner measure of the set $\{X\ne Y\bigr\}$ is 0. Hence, it suffices to let
$$\mathcal{A}=\sigma\bigl(\mathcal{C}\cup\bigl\{X\ne Y\bigr\}\bigr)$$
and to take $\mathbb{P}$ as the only extension of $\mathbb{Q}$ to $\mathcal{A}$ such that $\mathbb{P}(X\ne Y)=0$.

\medskip

\noindent Suppose now that $\mu|\mathcal{D}\ne\nu|\mathcal{D}$. Define $\lambda$, $f$ and $g$ by \eqref{f67y9k} and
\begin{gather*}
\gamma(A)=\frac{1}{\norm{\mu-\nu}_{\mathcal{D}}}\,\,\int_A\left(g-f\right)^+\,d\lambda\quad\quad\text{for all }A\in\mathcal{E}.
\end{gather*}
Since $\int\left(g-f\right)^+\,d\lambda=\norm{\mu-\nu}_\mathcal{D}$, such a $\gamma$ is a probability measure on $\mathcal{E}$. Let $(\Phi,\mathcal{C},\mathbb{Q})$ be any probability space which supports three independent random variables $U,\,X,\,Z$ with $U$ uniformly distributed on $(0,1)$ and
$$\mathbb{Q}(X\in A)=\mu(A)\quad\text{and}\quad\mathbb{Q}(Z\in A)=\gamma(A)\quad\text{for all }A\in\mathcal{E}.$$
Define
\begin{gather*}
G=\bigl\{f(X)\,U>g(X)\bigr\},\quad Y=Z\text{ on }G\quad\text{and}\quad Y=X\text{ on }G^c.
\end{gather*}
Then,
\begin{gather*}
\mathbb{Q}(G)=\mathbb{Q}\left[f(X)>g(X),\,U>\frac{g(X)}{f(X)}\right]=\int_{\{f>g\}}\left(1-\frac{g}{f}\right)\,f\,d\lambda
\\=\int_{\{f>g\}}(f-g)\,d\lambda=\norm{\mu-\nu}_{\mathcal{D}}.
\end{gather*}
Moreover, for each $A\in\mathcal{E}$,
\begin{gather*}
\mathbb{Q}(Y\in A)=\mathbb{Q}\bigl(G\cap\{Z\in A\}\bigr)+\mathbb{Q}\bigl(G^c\cap\{X\in A\}\bigr)
\medskip\\\medskip=\mathbb{Q}(G)\,\mathbb{Q}(Z\in A)+\mathbb{Q}\left[f(X)\,U\le g(X),\,X\in A\right]
\\=\int_A (g-f)^+d\lambda+\int_{A\cap\{f>g\}}\frac{g}{f}\,d\mu+\mu\bigl(A\cap\{f\le g\}\bigr).
\end{gather*}
If $A\in\mathcal{D}$, since $f=\frac{d\,(\mu|\mathcal{D})}{d\,(\lambda|\mathcal{D})}$, one obtains
\begin{gather*}
\int_{A\cap\{f>g\}}\frac{g}{f}\,d\mu+\mu\bigl(A\cap\{f\le g\}\bigr)=\int_{A\cap\{f>g\}}\frac{g}{f}\,f\,d\lambda+\int_{A\cap\{f\le g\}}f\,d\lambda=\int_A (f\wedge g)\,d\lambda.
\end{gather*}
Therefore,
\begin{gather*}
\mathbb{Q}(Y\in A)=\int_A (g-f)^+d\lambda+\int_A (f\wedge g)\,d\lambda=\int_A g\,d\lambda=\nu(A)\quad\quad\text{for each }A\in\mathcal{D}.
\end{gather*}
Finally, denoting by $\mathbb{Q}_*$ and $\mathbb{Q}^*$ the inner and outer measures corresponding to $\mathbb{Q}$, one obtains
$$\norm{\mu-\nu}_{\mathcal{D}}\le\mathbb{Q}_*(X\ne Y)\le\mathbb{Q}^*(X\ne Y)\le\mathbb{Q}(G)=\norm{\mu-\nu}_{\mathcal{D}}.$$
Therefore, to conclude the proof, it suffices to take $(\Phi,\mathcal{A},\mathbb{P})$ as the completion of $(\Phi,\mathcal{C},\mathbb{Q})$.
\end{proof}

\medskip

\section{Results}\label{res}

\noindent It is quite intuitive that, when investigating strong duality, the partition of $\Omega$ in the equivalence classes of $E$ plays a role. Let $\Pi$ denote such a partition, i.e.
$$\Pi=\bigl\{[x]:x\in\Omega\bigr\}\quad\quad\text{where }[x]=\bigl\{y\in\Omega:(x,y)\in E\bigr\}.$$
The $\sigma$-fields $\mathcal{G}_0$ and $\mathcal{G}_1$, introduced in Section \ref{intro}, can be written as
\begin{gather*}
\mathcal{G}_0=\bigl\{A\in\mathcal{F}:A\text{ is a union of elements of }\Pi\bigr\},
\\\quad\mathcal{G}_1=\bigl\{A\in\widehat{\mathcal{F}}:A\text{ is a union of elements of }\Pi\bigr\},
\end{gather*}
where $\widehat{\mathcal{F}}$ denotes the universally measurable $\sigma$-field with respect to $\mathcal{F}$. Another useful fact is
\begin{gather}\label{fuits33}
1_A(x)-1_A(y)\le 1-1_E(x,y)\quad\quad\text{for all }(x,y)\in\Omega\times\Omega
\end{gather}
provided the set $A\subset\Omega$ is a union of elements of $\Pi$.

\medskip

\noindent Our starting point is the following.

\begin{thm}\label{thm1} If $(\Omega,\mathcal{F})$ is a standard Borel space, then
\begin{gather*}
\inf_{P\in\Gamma(\mu,\nu)}(1-P(E))=\norm{\mu-\nu}_{\mathcal{G}_1}\quad\quad\text{for all }\mu,\,\nu\in\mathcal{P}(\mathcal{F}).
\end{gather*}
\end{thm}

\begin{proof} Let $\mu,\,\nu\in\mathcal{P}(\mathcal{F})$. In the notation of Section \ref{prel}, let $c=1-1_E$ and
$$L=\bigl\{(f,g):\,f,\,g\in M(\mathcal{F})\text{ and }f(x)+g(y)\le c(x,y)\text{ for all }(x,y)\in\Omega\times\Omega\bigr\}.$$
Since $(\Omega,\mathcal{F})$ is standard Borel, $\mu$ and $\nu$ are perfect. Therefore,
\begin{gather*}
\inf_{P\in\Gamma(\mu,\nu)}(1-P(E))=\inf_{P\in\Gamma(\mu,\nu)}P(c)=\sup_{(f,g)\in L}\bigl\{\mu(f)+\nu(g)\bigr\}.
\end{gather*}

\medskip

\noindent Given $(f,g)\in L$, define
$$\phi=\bigl(f-\sup f+1)^+\quad\text{and}\quad\psi=g+\sup f-1.$$
On noting that
$$\sup f+\sup g=\sup_{(x,y)}\{f(x)+g(y)\}\le\sup c\le 1,$$
one obtains $(\phi,\psi)\in L$. Moreover, $0\le\phi\le 1$ and $\mu(\phi)+\nu(\psi)\ge\mu(f)+\nu(g)$. Hence,
$$\inf_{P\in\Gamma(\mu,\nu)}(1-P(E))=\sup_{(f,g)\in L}\bigl\{\mu(f)+\nu(g)\bigr\}=\sup_{
                                                                               \begin{array}{c}
                                                                                 (f,g)\in L \\
                                                                                 0\le f\le 1 \\
                                                                               \end{array}}\,\{\mu(f)+\nu(g)\}.$$
\medskip

\noindent Next, fix $\epsilon>0$ and take $(f,g)\in L$ such that $0\le f\le 1$ and
$$\mu(f)+\nu(g)+\epsilon > \inf_{P\in\Gamma(\mu,\nu)}(1-P(E)).$$
Define
$$h(x)=\sup_{y\in [x]}f(y)$$
and note that $h(x)+g(y)\le c(x,y)$ for all $(x,y)$. Letting $y=x$, one obtains
$$g(x)\le c(x,x)-h(x)=-h(x)\quad\quad\text{for all }x\in\Omega.$$
Since $h(x)=h(y)$ whenever $(x,y)\in E$, for each $a\in\mathbb{R}$ the set $\bigl\{h>a\bigr\}$ is a union of elements of $\Pi$. Moreover, since $(\Omega,\mathcal{F})$ is standard Borel, the projection theorem yields
$$\bigl\{h>a\bigr\}=\bigl\{x\in\Omega:(x,y)\in E\text{ and }f(y)>a\text{ for some }y\in\Omega\bigr\}\in\widehat{\mathcal{F}};$$
see e.g. Theorem A1.4, page 562, of \cite{KAL}. Hence, $\bigl\{h>a\bigr\}\in\mathcal{G}_1$. To sum up,
$$h\in M(\mathcal{G}_1),\quad 0\le h\le 1,\quad h\ge f,\quad -h\ge g.$$
Therefore,
\begin{gather*}
\norm{\mu-\nu}_{\mathcal{G}_1}=\sup_{
                                                                               \begin{array}{c}
                                                                                 f\in M(\mathcal{G}_1) \\
                                                                                 0\le f\le 1 \\
                                                                               \end{array}}\,\abs{\mu(f)-\nu(f)}\ge\mu(h)-\nu(h)
\\\ge\mu(f)+\nu(g) > \inf_{P\in\Gamma(\mu,\nu)}(1-P(E))-\epsilon.
\end{gather*}

\medskip

\noindent Finally, fix $P\in\Gamma(\mu,\nu)$ and $A\in\mathcal{G}_1$. Since $A$ is a union of elements of $\Pi$, inequality \eqref{fuits33} yields
$$1-P(E)\ge\int\bigl\{1_A(x)-1_A(y)\bigr\}\,P(dx,dy)=\mu(A)-\nu(A).$$
Hence,
$$\inf_{P\in\Gamma(\mu,\nu)}(1-P(E))\ge\norm{\mu-\nu}_{\mathcal{G}_1}$$
and this concludes the proof.

\end{proof}

\medskip

\noindent If regarded as a tool to get strong duality, Theorem \ref{thm1} has two gaps:

\medskip

\begin{itemize}

\item $\mathcal{G}_1$ is not a sub-$\sigma$-field of $\mathcal{F}$;

\item Theorem \ref{thm1} involves the $\inf$ and not the $\min$ over $\Gamma(\mu,\nu)$.

\end{itemize}

\medskip

\noindent The rest of this note focus on these two points.

\bigskip

\subsection{$\mathcal{G}_1$ is not a sub-$\sigma$-field of $\mathcal{F}$.}\label{v66t8}

\noindent The fact that $\mathcal{G}_1$ is not a sub-$\sigma$-field of $\mathcal{F}$ may be seen as unsuitable. Hence, we now prove the following result.

\begin{thm}\label{ncp97u}
If $(\Omega,\mathcal{F})$ is a standard Borel space, then, for all $\mu,\,\nu\in\mathcal{P}(\mathcal{F})$, there is a set $B\in\mathcal{F}$ such that
\begin{gather*}
\mu(B)=\nu(B)=1\quad\text{and}\quad\inf_{P\in\Gamma(\mu,\nu)}(1-P(E))=\norm{\mu-\nu}_{\mathcal{G}_B}
\end{gather*}
where $\mathcal{G}_B=\bigl\{A\in\mathcal{F}:1_A(x)=1_A(y)\text{ for all }(x,y)\in E\cap (B\times B)\bigr\}$.
\end{thm}

\begin{proof}
Let $\mu,\,\nu\in\mathcal{P}(\mathcal{F})$. By \eqref{g8m4e3} and Theorem \ref{thm1}, there is $D\in\mathcal{G}_1$ such that
\begin{gather*}
\inf_{P\in\Gamma(\mu,\nu)}(1-P(E))=\norm{\mu-\nu}_{\mathcal{G}_1}=\mu(D)-\nu(D).
\end{gather*}
Since $D$ is universally measurable with respect to $\mathcal{F}$, there is $A\in\mathcal{F}$ such that
$$\frac{\mu+\nu}{2}\bigl(A\Delta D\bigr)=0,$$
or equivalently $\mu(A\Delta D)=\nu(A\Delta D)=0$. Let
$$T=\bigl\{(x,y)\in E:1_A(x)\neq 1_A(y)\bigl\}.$$
Since $D$ is a union of elements of $\Pi$, then $1_D(x)=1_D(y)$ for all $(x,y)\in E$. Hence,
$$P(T)=P\bigl\{(x,y)\in E:1_D(x)\neq 1_D(y)\bigl\}=P(\emptyset)=0\quad\text{for each }P\in\Gamma(\mu,\nu)$$
where the first equality is because $\mu(A\Delta D)=\nu(A\Delta D)=0$. Since $(\Omega,\mathcal{F})$ is standard Borel and $P(T)=0$ for all $P\in\Gamma(\mu,\nu)$, there is $B\in\mathcal{F}$ such that $\mu(B)=\nu(B)=1$ and
$$T\subset (B^c\times\Omega)\cup (\Omega\times B^c);$$
see \cite{HS} and \cite[p. 2345]{RIG20}. Therefore $A\in\mathcal{G}_B$, which in turn implies
$$\inf_{P\in\Gamma(\mu,\nu)}(1-P(E))=\mu(D)-\nu(D)=\mu(A)-\nu(A)\le\norm{\mu-\nu}_{\mathcal{G}_B}.$$
To prove the reverse inequality, fix any $C\in\mathcal{G}_B$ and $P\in\Gamma(\mu,\nu)$. Then,
$$P(B\times B)=1\quad\text{and}\quad 1_C(x)-1_C(y)\le 1-1_E(x,y)\quad\text{for all }(x,y)\in B\times B.$$
Hence,
$$\mu(C)-\nu(C)=\int \bigl(1_C(x)-1_C(y)\bigr)\,P(dx,dy)\le 1-P(E),$$
which in turn implies $\norm{\mu-\nu}_{\mathcal{G}_B}\le\inf_{P\in\Gamma(\mu,\nu)}(1-P(E))$.
\end{proof}

\medskip

\noindent As already noted, the advantage of Theorem \ref{ncp97u} with respect to Theorem \ref{thm1} is that $\mathcal{G}_B$ is a sub-$\sigma$-field of $\mathcal{F}$ while $\mathcal{G}_1$ is not. The disadvantage is that $\mathcal{G}_B$ is not universal, for it depends on the pair $(\mu,\nu)$. However, for fixed $(\mu,\nu)$, since $\mu(B)=\nu(B)=1$ and $E\cap (B\times B)$ is a measurable equivalence relation on $B$, it is reasonable to replace $\Omega$ with $B$ and $E$ with $E\cap (B\times B)$. In other terms, for fixed $(\mu,\nu)$, it makes sense to involve $\mathcal{G}_B$ in the notion of strong duality.

\bigskip

\subsection{Existence of primal minimizers}

\noindent Quite surprisingly, in mass transportation theory, existence of primal minimizers seems to have received only a little attention to date; see e.g. \cite{BLS} and \cite{JAFFE}. To our knowledge, when the cost $c$ is not lower semi-continuous, the only available results are in \cite{KELL} and require $c$ to be suitably approximable by regular costs. However, such results do not apply to our case where $c=1-1_E$.

\medskip

\noindent Let $(\Omega,\mathcal{F})$ be a standard Borel space and $c=1-1_E$. Then, $c$ is lower semi-continuous if and only if $E$ is closed, and in this case $E$ is strongly dualizable. Similarly, $E$ is strongly dualizable if the elements of the partition $\Pi$ are the atoms of a countably generated sub-$\sigma$-field of $\mathcal{F}$; see \cite{JAFFE} and \cite{PR2023}. As noted above, however, we are not aware of any (reasonable) condition for a primal minimizer to exist. In the sequel, we discuss two strategies for circumventing this problem.

\medskip

\noindent The first strategy is possibly expected and lies in using finitely additive probabilities. Let
$$M(\mu,\nu)=\bigl\{\text{finitely additive probabilities on }\mathcal{F}\otimes\mathcal{F}\text{ with marginals }\mu\text{ and }\nu\bigr\}.$$

\begin{thm}\label{f59ik1}
Let $(\Omega,\mathcal{F})$ be a standard Borel space. Then,
\begin{gather*}
\min_{P\in M(\mu,\nu)}(1-P(E))=\norm{\mu-\nu}_{\mathcal{G}_1}\quad\quad\text{for all }\mu,\,\nu\in\mathcal{P}(\mathcal{F}).
\end{gather*}
Moreover, for all $\mu,\,\nu\in\mathcal{P}(\mathcal{F})$ there is $B\in\mathcal{F}$ such that
\begin{gather*}
\mu(B)=\nu(B)=1\quad\text{and}\quad\min_{P\in M(\mu,\nu)}(1-P(E))=\norm{\mu-\nu}_{\mathcal{G}_B}.
\end{gather*}
\end{thm}

\begin{proof}
Just apply Theorems \ref{thm1} and \ref{ncp97u} and note that, by Theorem 2 of \cite{RIG23},
$$\min_{P\in M(\mu,\nu)}(1-P(E))=\inf_{P\in\Gamma(\mu,\nu)}(1-P(E)).$$
\end{proof}

\medskip

\noindent A remark on Theorem \ref{f59ik1} is in order. Let $P$ be a finitely additive primal minimizer, in the sense that $P\in M(\mu,\nu)$ and $1-P(E)=\inf_{Q\in\Gamma(\mu,\nu)}(1-Q(E))$. Moreover, let $\mathcal{R}$ be the field generated by the measurable rectangles $A\times B$ with $A,\,B\in\mathcal{F}$. Since $\mu$ and $\nu$ are perfect (due to $(\Omega,\mathcal{F})$ is standard Borel), the restriction $P|\mathcal{R}$ is $\sigma$-additive; see e.g. \cite{R1996}. Hence, it is tempting to define $P'$ as the only $\sigma$-additive extension of $P|\mathcal{R}$ to $\sigma(\mathcal{R})=\mathcal{F}\otimes\mathcal{F}$. Then, $P'\in\Gamma(\mu,\nu)$ but it is {\em not} necessarily true that $P'(E)=P(E)$. Hence, $P'$ needs not be a primal minimizer.

\medskip

\noindent The second strategy for dealing with primal minimizers is summarized by the next result.

\begin{thm}\label{e55t7b8} For all $\mu,\,\nu\in\mathcal{P}(\mathcal{F})$, there is $\nu_0\in\mathcal{P}(\mathcal{F})$ such that
$$\nu_0=\nu\text{ on }\mathcal{G}_0\quad\text{and}\quad\min_{P\in\Gamma(\mu,\nu_0)}(1-P(E))=\norm{\mu-\nu}_{\mathcal{G}_0}$$
where $\mathcal{G}_0$ is the $\sigma$-field introduced in Section \ref{intro}.
\end{thm}

\begin{proof} By Lemma \ref{kl8}, applied with $S=\Omega$, $\mathcal{E}=\mathcal{F}$ and $\mathcal{D}=\mathcal{G}_0$, there are a probability space $(\Phi,\mathcal{A},\mathbb{P})$ and two measurable maps $X,\,Y:(\Phi,\mathcal{A})\rightarrow (\Omega,\mathcal{F})$ such that
\begin{gather*}
P(X\in A)=\mu(A)\text{ for all }A\in\mathcal{F},\quad P(Y\in A)=\nu(A)\text{ for all }A\in\mathcal{G}_0,
\\\bigl\{X\ne Y\bigr\}\in\mathcal{A}\quad\text{and}\quad\mathbb{P}(X\ne Y)=\norm{\mu-\nu}_{\mathcal{G}_0}.
\end{gather*}
Up to replacing $(\Phi,\mathcal{A},\mathbb{P})$ with its completion, it can be assumed that $(\Phi,\mathcal{A},\mathbb{P})$ is complete. Let $\mathbb{P}_*$ and $\mathbb{P}^*$ be the inner and outer measures corresponding to $\mathbb{P}$. Because of \eqref{fuits33},
$$1_{\{X\in A\}}-1_{\{Y\in A\}}\le 1_{\{(X,Y)\notin E\}}\le 1_{\{X\neq Y\}}\quad\quad\text{for each }A\in\mathcal{G}_0.$$
Therefore,
\begin{gather*}
\mu(A)-\nu(A)=\int\bigl(1_{\{X\in A\}}-1_{\{Y\in A\}}\bigr)\,d\mathbb{P}\le\mathbb{P}_*\bigl((X,Y)\notin E\bigr)
\\\le\mathbb{P}^*\bigl((X,Y)\notin E\bigr)\le\mathbb{P}(X\neq Y)=\norm{\mu-\nu}_{\mathcal{G}_0}
\end{gather*}
for each $A\in\mathcal{G}_0$, which in turn implies
$$\norm{\mu-\nu}_{\mathcal{G}_0}\le\mathbb{P}_*\bigl((X,Y)\notin E\bigr)\le\mathbb{P}^*\bigl((X,Y)\notin E\bigr)\le\norm{\mu-\nu}_{\mathcal{G}_0}.$$
Since $(\Phi,\mathcal{A},\mathbb{P})$ is complete, one obtains
$$\bigl\{(X,Y)\notin E\bigr\}\in\mathcal{A}\quad\text{and}\quad\mathbb{P}\bigl((X,Y)\notin E\bigr)=\norm{\mu-\nu}_{\mathcal{G}_0}.$$
To conclude the proof, note that $\bigl\{(X,Y)\in H\bigr\}\in\mathcal{A}$ for each $H\in\mathcal{F}\otimes\mathcal{F}$ and define
$$\nu_0(A)=\mathbb{P}(Y\in A)\quad\text{and}\quad P(H)=\mathbb{P}\bigl((X,Y)\in H\bigr)\quad\quad\text{for all }A\in\mathcal{F}\text{ and }H\in\mathcal{F}\otimes\mathcal{F}.$$
Then, $\nu_0=\nu$ on $\mathcal{G}_0$, $P\in\Gamma(\mu,\nu_0)$ and
$$1-P(E)=\norm{\mu-\nu}_{\mathcal{G}_0}\le 1-Q(E)\quad\quad\text{for each }Q\in\Gamma(\mu,\nu_0).$$
\end{proof}

\medskip

\noindent It is worth noting that, in Theorem \ref{e55t7b8}, $(\Omega,\mathcal{F})$ is not required to be a standard Borel space. In addition, Theorem \ref{e55t7b8} has the following useful consequence.

\begin{cor}\label{vr44e8u1}
Let $(\Omega,\mathcal{F})$ be a standard Borel space. If $E\in\widehat{\mathcal{F}}\otimes\mathcal{G}_1$, then
$$\min_{P\in\Gamma(\mu,\nu)}\,(1-P(E))=\norm{\mu-\nu}_{\mathcal{G}_1}\quad\quad\text{for all }\mu,\,\nu\in\mathcal{P}(\mathcal{F}).$$
Moreover, if $E\in\mathcal{F}\otimes\mathcal{G}_0$, one also obtains
$$\min_{P\in\Gamma(\mu,\nu)}\,(1-P(E))=\norm{\mu-\nu}_{\mathcal{G}_0}\quad\quad\text{for all }\mu,\,\nu\in\mathcal{P}(\mathcal{F}).$$
\end{cor}

\begin{proof}
We first recall that, for each $\gamma\in\mathcal{P}(\mathcal{F})$, the only extension of $\gamma$ to $\widehat{\mathcal{F}}$ is still denoted by $\gamma$. Moreover, since $(\Omega,\mathcal{F})$ is standard Borel, every probability measure on $\widehat{\mathcal{F}}$ is perfect.

\medskip

\noindent We begin with the second part of the corollary, so that we assume $E\in\mathcal{F}\otimes\mathcal{G}_0$. Given $\mu,\,\nu\in\mathcal{P}(\mathcal{F})$, by Theorem \ref{e55t7b8}, there are $\nu_0\in\mathcal{P}(\mathcal{F})$ and $P_0\in\Gamma(\mu,\nu_0)$ such that $\nu_0=\nu$ on $\mathcal{G}_0$ and $1-P_0(E)=\norm{\mu-\nu}_{\mathcal{G}_0}$. In addition, since $\mu$ and $\nu$ are perfect, by Theorem 9 of \cite{R1996}, there is $P\in\Gamma(\mu,\nu)$ such that $P=P_0$ on $\mathcal{F}\otimes\mathcal{G}_0$. Since $E\in\mathcal{F}\otimes\mathcal{G}_0$, one obtains $P(E)=P_0(E)$. Therefore,
$$P\in\Gamma(\mu,\nu)\quad\text{and}\quad 1-P(E)=\norm{\mu-\nu}_{\mathcal{G}_0}\le 1-Q(E)\quad\quad\text{for each }Q\in\Gamma(\mu,\nu)$$
where the inequality is by \eqref{fuits33}.

\medskip

\noindent Next, we assume $E\in\widehat{\mathcal{F}}\otimes\mathcal{G}_1$. Let $\widehat{\Gamma}(\mu,\nu)$ be the collection of probability measures $\widehat{P}$ on $\widehat{\mathcal{F}}\otimes\widehat{\mathcal{F}}$ such that
$$\widehat{P}(A\times\Omega)=\mu(A)\quad\text{and}\quad\widehat{P}(\Omega\times A)=\nu(A)\quad\text{for all }A\in\widehat{\mathcal{F}}.$$
Since $E\in\widehat{\mathcal{F}}\otimes\mathcal{G}_1$ and $\mu$ and $\nu$ are perfect (where $\mu$ and $\nu$ are now regarded as probability measures on $\widehat{\mathcal{F}}$) the first part of this proof can be repeated with $(\Omega,\widehat{\mathcal{F}})$ and $\mathcal{G}_1$ in the place of $(\Omega,\mathcal{F})$ and $\mathcal{G}_0$. Hence,
$$1-\widehat{P}(E)=\norm{\mu-\nu}_{\mathcal{G}_1}\quad\text{for some}\quad\widehat{P}\in\widehat{\Gamma}(\mu,\nu).$$
Finally, denoting by $P$ the restriction of $\widehat{P}$ on $\mathcal{F}\otimes\mathcal{F}$, one obtains
$$P\in\Gamma(\mu,\nu)\quad\text{and}\quad 1-P(E)=\norm{\mu-\nu}_{\mathcal{G}_1}\le 1-Q(E)\quad\quad\text{for each }Q\in\Gamma(\mu,\nu).$$
\end{proof}

\medskip

\noindent When $E\in\mathcal{F}\otimes\mathcal{G}_0$, Corollary \ref{vr44e8u1} slightly improves \cite[Theo. 3.13]{JAFFE} which requires $E\in\mathcal{G}_0\otimes\mathcal{G}_0$. Instead, when $E\in\widehat{\mathcal{F}}\otimes\mathcal{G}_1$, we are not aware of any analogous result. In any case, here is an application of Corollary \ref{vr44e8u1}.

\medskip

\begin{ex}
Let $\Omega$ be a Polish space and $\mathcal{F}=\mathcal{B}(\Omega)$. A subset of $\Omega$ is a $G_\delta$ if it is a countable intersection of open sets. In particular, open and closed subsets of $\Omega$ are both $G_\delta$. Moreover, by Corollary \ref{vr44e8u1}, {\em if the equivalence classes of $E$ are $G_\delta$, then}
\begin{gather*}
\min_{P\in\Gamma(\mu,\nu)}\,(1-P(E))=\norm{\mu-\nu}_{\mathcal{G}_1}\quad\quad\text{for all }\mu,\,\nu\in\mathcal{P}(\mathcal{F}).
\end{gather*}
\noindent To prove the latter claim, it suffices to show that $E\in\widehat{\mathcal{F}}\otimes\mathcal{G}_1$. For $A\in\mathcal{F}$, define
$$A^*=\{x\in\Omega:\,\exists\, y\in A\text{ such that }(x,y)\in E\}.$$
Since $A^*$ is the projection on the $x$-axis of the Borel set $E\cap (\Omega\times A)$, the projection theorem yields $A^*\in\widehat{\mathcal{F}}$; see again Theorem A1.4, page 562, of \cite{KAL}. Since $A^*$ is a union of equivalence classes of $E$, one also obtains $A^*\in\mathcal{G}_1$. Having noted this fact, fix a countable basis $\mathcal{U}$ for the topology of $\Omega$ and define
$$\mathcal{V}=\sigma(A^*:A\in\mathcal{U}).$$
Then, $\mathcal{V}$ is countably generated and $\mathcal{V}\subset\mathcal{G}_1$. If $A$ and $B$ are any disjoint $G_\delta$ sets, there is $U\in\mathcal{U}$ such that
$$A\cap U\ne\emptyset\text{ and }B\cap U=\emptyset\quad\text{or}\quad A\cap U=\emptyset\text{ and }B\cap U\ne\emptyset;$$
see the proof of Lemma 2 in \cite{MILLER}. Hence, if $A$ and $B$ are two disjoint equivalence classes of $E$, then
$$A\subset U^*\text{ and }B\cap U^*=\emptyset\quad\text{or}\quad A\cap U^*=\emptyset\text{ and }B\subset U^*$$
for some $U\in\mathcal{U}$. This implies that the equivalence classes of $E$ are precisely the atoms of $\mathcal{V}$. Finally, since $\mathcal{V}$ is countably generated, there is a function $f:\Omega\rightarrow\mathbb{R}$ such that $\mathcal{V}=\sigma(f)$. Therefore,
$$E=\bigl\{(x,y):f(x)=f(y)\bigr\}\in\mathcal{V}\otimes\mathcal{V}\subset\widehat{\mathcal{F}}\otimes\mathcal{G}_1.$$
\end{ex}

\medskip

\noindent We close this note with a last result. While not practically useful, it still provides some information on primal minimizers.

\begin{thm}\label{njil9}
Let $(\Omega,\mathcal{F})$ be a standard Borel space and $P\in\Gamma(\mu,\nu)$ where $\mu,\,\nu\in\mathcal{P}(\mathcal{F})$. Then, $P$ is a primal minimizer (with respect to $c=1-1_E$) if and only if
\begin{gather}\label{cns4o}
P(E)=1-P(A\times A^c)\quad\quad\text{for some }A\in\mathcal{G}_1.
\end{gather}
\end{thm}

\begin{proof}
By \eqref{g8m4e3}, there is $A\in\mathcal{G}_1$ such that $\norm{\mu-\nu}_{\mathcal{G}_1}=\mu(A)-\nu(A)$. Hence, if $P$ is a primal minimizer, Theorem \ref{thm1} implies
\begin{gather*}
1-P(E)=\norm{\mu-\nu}_{\mathcal{G}_1}=\mu(A)-\nu(A)=\int\bigl(1_A(x)-1_A(y)\bigr)\,P(dx,dy)
\\\le P(A\times A^c)\le P(E^c)=1-P(E).
\end{gather*}
Hence, condition \eqref{cns4o} holds. Conversely, if \eqref{cns4o} holds for some $A\in\mathcal{G}_1$, then
\begin{gather*}
P\bigl\{(x,y):1_A(x)-1_A(y)=1-1_E(x,y)\bigr\}=P(E)+P\bigl\{(x,y)\in E^c:1_A(x)-1_A(y)=1\bigr\}
\\=P(E)+P\bigl(E^c\cap (A\times A^c)\bigr)=P(E)+P(A\times A^c)=1.
\end{gather*}
Therefore, for each $Q\in\Gamma(\mu,\nu)$,
\begin{gather*}
1-P(E)=\int\bigl(1_A(x)-1_A(y)\bigr)\,P(dx,dy)=\mu(A)-\nu(A)
\\=\int\bigl(1_A(x)-1_A(y)\bigr)\,Q(dx,dy)\le 1-Q(E)
\end{gather*}
where the last inequality is by \eqref{fuits33}. Hence, $P$ is a primal minimizer.
\end{proof}

\end{document}